\documentclass[12pt]{article}
\usepackage{amsmath,amssymb,amsfonts,euscript}
\usepackage{graphicx,hyperref,refcheck}
\newtheorem{thm}{\bf Theorem}
\newtheorem{lem}[thm]{\bf Lemma}
\newtheorem{cor}[thm]{\bf Corollary}

\newtheorem{dfn}[thm]{\bf Definition}
\newcommand{\dok}{{\bf Proof.\ }}

\newcommand{\rr}{\mathbb{R}}

\newcommand {\al} {\alpha}
\newcommand {\be} {\beta}
\newcommand {\da} {\delta}

\newcommand {\ga} {\gamma}

\newcommand {\la} {\lambda}
\newcommand {\La} {\Lambda}

\newcommand {\sa} {\sigma}

\newcommand{\e}{\varepsilon}

\newcommand{\IN}{{\subset}}

\newcommand{\mC}{{\mathsf C}}
\newcommand \w {\widetilde}
\newcommand \yy {^{-1}}

\newcommand \dd  {\partial}
\newcommand {\mmm}{{\setminus}}

\newcommand{\8}{{\infty}}

\newcommand{\ia}{{I^*}}
\newcommand{\0}{{\varnothing}}
\newcommand{\vse}{$\blacksquare$}

\newcommand{\bj}{{\bf {j}}}
\newcommand{\bi}{{\bf {i}}}

\newcommand{\eS}{{\EuScript S}}

\newcommand{\eP}{{\EuScript P}}

\newcommand{\eG}{{\EuScript G}}

\newcommand{\eF}{{\EuScript F}}
\newcommand{\eZ}{{\EuScript Z}}

\def \fix {\mathop{\rm fix}\nolimits}
\def \Lip {\mathop{\rm Lip}\nolimits}

\def \Id {\mathop{\rm Id}\nolimits}

\title{On weak separation property for self-affine Jordan arcs.}

\author{Olesya Chelkanova \and 
Andrey Tetenov}

\begin{document}

\maketitle

\newcommand{\beq}{\begin{equation}}
\newcommand{\eeq}{\end{equation}}
\newcommand{\bpm}{\begin{pmatrix}}
\newcommand{\epm}{\end{pmatrix}}

 \begin{abstract} 
 We consider self-affine arcs in $\mathbb R^2$ and prove that violation of "inner" weak separation property for such arcs implies that the arc is a parabolic segment. Therefore, if a self-affine Jordan arc is not a parabolic segment, then it is the attractor of some multizipper.
 \end{abstract}


\section{Introduction}

The idea of associated family of similarities for  a system $\eS=\{ S_1,...,S_m \}$ of similarities in $ \rr^d $ was initially proposed by C.Bandt and S.Graf \cite{SSS7} to analyse the measure and dimension properties of the attractor $ K $  of the system $ \eS $. This approach was developed in \cite{Schief} to result in Weak Separation Condition \cite{Edgdas,Lau,DE, Zer}. Violation of WSC results not only in the measure drop for $K$ \cite{KT2F, TKV} in its dimension, but it also implies some special geometric properties of $K$  and  rigidity phenomena for the deformations of self-similar structure on $K$\cite{BR,Trg,TCh}.

Though this scope of ideas and methods initially had self-similar sets as its target, there always was an attractive idea to extend it to more general classes of self-similar sets.

We consider how Weak Separation Condition (or its violation) applies to self-affine Jordan arcs in plane and show that  structure and rigidity  theorems for self-similar Jordan arcs \cite{ATK,Atet1} have their self-affine analogues.

The main result  of the current paper is the following

\begin{thm}\label{main}
Let $\ga$ be a self-affine Jordan arc in $\rr^2$ whish is not a parabolic segment. Then  $\ga$ is a component of the attractor of some self-affine multizipper $\eZ$ 
  \end{thm}
  
 As a main step for this result we prove the following rigidity theorem for a very general class of self-affine arcs, which need not be finitely generated:
  
  \begin{thm}\label{gen} 
 
 Let   $\ga=\ga(a_0,a_1)$ be a Jordan arc with endpoints $a_0, a_1$ in $\rr^2$ such that\\
 (i) For any $\e>0$ and for any non-degenerate subarc $\ga' \IN \ga$ there is an affine map $S$ such that $S (\ga) \IN \ga'$ and $\Lip S < \e$\\ 
 (ii) There is a sequence of affine maps $f_k$ converging to $\Id$ such that $f_k(\ga) \cap \ga = \ga(f(a_0),a_1)$ and $\fix(f_k)\cap\ga=\0$;\\ Then $\ga$ is a parabolic segment.
  \end{thm}
  
  In finitely generated case this theorem becomes

\begin{thm}\label{fin}
Let a Jordan arc $\ga\IN\rr^2$ with endpoints  $a_0,a_1$
be the attractor of a system  $ \eS =\{S_1, ..., S_m\}$ of contracting
affine maps. Let $ {\eF(\eS)}$ be the associated family for the system $\eS$. If there is a sequence
 $f_n \in {\eF(\eS)} \mmm  \{ Id \}$ such that $f_n \to Id$, and $f_n(\ga)\cap\ga\neq\0$  then $\ga$ is a parabolic segment.
\end{thm}

The proof of Theorems  2 and 3 uses the result of C.~Bandt and A.~S.~Kravchenko  \cite{BK} that
except for parabolic arcs and segments, there are no twice continuously differentiable
self-affine curves in the plane.\\

{\bf 1.Definitions and notation.}\\

Let $\eS=\{S_1,\dots,S_m\}$ be a system of contracting affine maps in $\mathbb R^d$. The unique nonempty compact set $K=K(\eS)$ such that $K=\bigcup\limits_{i=1}^{m} S_i (K)$, is called the {\em attractor} of the system $\eS$, or a {\em self-affine set} generated by the system $\eS$.

A system $\eS$  is {\em irreducible} if, for every proper subsystem 
$\eS'\IN \eS$,
the attractor of $\eS'$ is different from
the attractor of the system $\eS$.


By $I=\{1,2,...,m\}$ we denote the set of indices, $\ia=\bigcup\limits_{n=1}^\8 I^n$  
is the set of all  multiindices $\bi=i_1i_2...i_n$, and we denote $S_\bi = S_{i_1} S_{i_2} ... S_{i_n}$. 
The set of all infinite sequences $I^{\8}=\{{\bf \al}=\al_1\al_2\ldots,\ \ \al_i\in I\}$ is  the
{\em index space}; and $\pi:I^{\8}\rightarrow K$ is the {\em index map}, which maps a sequence $\bf\al$ to  the point $\bigcap\limits_{n=1}^\8 K_{\al_1\ldots\al_n}$.

The set $\eF$ of all compositions $S_{\bf j}^{-1}S_{\bf i}$, where ${\bf i}, {\bf j} \in I^{*}$ and $i_1\neq j_1$ is called the associated family of affine mappings for the system $\eS$. The system $\eS$ has the {\em weak separation property} (WSP) if and only if $\rm Id \notin \overline{\eF \setminus \rm Id}$. \\

If $\ga$ is a Jordan arc with endpoints $a_0,a_1$, we denote  its subarc $\ga'$ with endpoints $x,y\in\ga$ by $\ga(x,y)$. We order the points in $\ga$ putting $a_0<a_1$ and write $x<y$ if $y\in \ga(x,a_1)$. We denote the diameter of a set $A$ by $|A|$.\\

 {\bf 2.  Representing $\ga$ as a limit of $\e$-nets $P(k,x)$.}\\

  Applying if necessary a coordinate change, we may suppose that the arc  $\ga$ lies in the unit  disc $D= \{x^2 +y^2 \le 1\}$.

  It follows from the condition (ii) that the subarcs $\sa_{k,0}= \ga \setminus f_k(\ga)$ and  $\sa_{k,1}=f_k(\ga) \setminus {f_k}^2(\ga)$ are disjoint. Proceeding by induction we get a sequence of subarcs 
   \beq \sa_{k,n} = f_k^n (\sa_{k,0})=f_k^n (\ga) \setminus f_k^{n+1}(\ga) \eeq
   which have endpoints  $f_k^n (a_0), f_k^{n+1}(a_0)$ and have disjont interiors as long as respective subarcs lie in $\ga$.  Since $f_k$ has no fixed points in $\ga$,  there is a maximal number $N_k$ for which  $\bigcup\limits_{n=0}^{N_k-1}\sa_{k,n}=\ga(a_0,f^{N_k}_k (a_0))\IN \ga$. 
   Let $\sa_{k,N_k} = f_k^{N_k} (\sa_{k,0}) \cap \ga=\ga(f^{N_k}_k (a_0),a_1)$.\\
  
 By the compactness of the arc  $\ga$ for any $\e> 0$ there is such $\da$, that if $x_1,x_2 \in \ga$ and $d(x_1,x_2 ) <\da$, then the diameter of the subarc $\ga(x_1,x_2)$ is less than $\e$.\\ 
  
  Therefore for any  $\e> 0$ there is such $N$, that if $k<N$ then $\|f_k (x) - x\|<\da$ for any $x\in \ga$, therefore the diameters of the subarcs $\sa_{k,n}$ are not greater than $\e$. \\
  
  For any $k$ and for any $x\in \ga$ the point $x$ lies in one of subarcs $f_k^{n_k}(\sa_{x,0})$. Denote $P(k,x)=\{f_k^n (x), -n_k \le n \le N_k-n_0\} $. Then Hausdorff distance between  $P(k,x)$ and $\ga$ is not greater than $\max\{|\sa_{k,n}|, 0\le n\le N_k\}$.\\ 

Therefore for any choice of the sequence $x_k \in \ga$  the sequence of sets $P(k, x_k)$ converges to $\ga$ in Hausdorff metrics.\\
  
  {\bf 3. Five types of affine maps and their associated vector fields.}
  
  Since the sequence $f_k$ converges to $\Id$, we suppose that all $f_k $ are sufficiently close to $\Id$ so that for any  $f_k$ we can correctly define its power $f_k^t, t\in \rr$, satisfying the conditions: \\
  
   1. For any $t_1,t_2\in\rr$, $f_k^{t_1} \circ f_k^{t_2}=f_k^{t_1+t_2}$;\qquad
  2. $f_k^0=\Id$ and  $f_k^1={f_k}$.\\

For that reason we divide the set of  non-degenerate affine maps $f (x) = Ax+b$ on $\rr^2$, where $A$ is a non-degenerate matrix and $b$ is a vector to five following types, depending on the eigenvalues $\la_1$ and $\la_2$ of the matrix $A$ and on the translation vector $b$:\\

{\bf Type 1.} If both eigenvalues $\la_1$ and $\la_2$ are not equal to $1$, then the map $f (x)$ has unique fixed point $x_0= {(E-A)}^{-1} b$.
By our assumptions,  $\|A-E\|<1$, therefore  $A=e^B$, where $B= \sum \limits_ {n=1}^\infty (-1)^{n+1} \dfrac{(A-E)^n}{n}$ is the matrix logarithm of $A$.
Since $f(x)= A(x-x_0)+x_0$, we put \beq \label{bfor1} f^t(x) = e^{B t}(x-x_0)+x_0\eeq

In this case for  any $x\neq x_0$, $\{ f^t(x), t\in \rr\}$ is an integral curve of  autonomous system $\dot x= B(x-x_0)$.\\

{\bf Types 2 and 3.} If  $\la_1 \neq 1$ and $\la_2 =1$ and $e_1,e_2$ are respective eigenvectors, then the map $f$ can be represented by $f(x)=Ax+ae_1+be_2$.

In this case the matrix logarithm $ B $ has eigenvalues $ \log\la_1 $ and $0$ and the equation \beq    f^t(x ) = e^{Bt} x + a \dfrac{\la_1^t -1}{\la_1-1}e_1 +b t e_2 
 \eeq defines some integral curve of the autonomous system
  \beq \label{dyn2}  \dot x = Bx + \dfrac{a \log\la_1}{\la_1-1}e_1 +b e_2  
 \eeq
We refer $f$ to the {\bf Type 2} if $b=0$. In this case
the right side in \eqref{dyn2} is a multiple of $e_1$, and integral curves are straight lines parallel to $e_1$.
If $ x =\dfrac{a }{1-\la_1}e_1+te_2$, then $Bx=- \dfrac{a \log\la_1}{\la_1-1}e_1$, so the right side in \eqref{dyn2} vanishes, and $L=\left\{\dfrac{a e_1}{1-\la_1} +e_2t, t\in\rr\right\}$
is the line consisting of fixed points of $f$.\\

$S$ is referred to {\bf Type 3} if $b\neq 0$. The system  \eqref{dyn2} has no fixed points in this case. The right side of \eqref{dyn2} on the line $L$ is equal to $be_2$, so $L$ is  the  invariant straight line.  The vector field is invariant under translations by $te_2, t\in\rr$, and there is the  minimal value for $ \| \dot x\| $ which is equal to $ |b| \|e_1\| |\sin \al_{12}| $, where
 $ \al_{12} $ is the angle between $e_1$ and $e_2$.\\

{\bf Type 4.} It is the case when the eigenvalues of $A$ are $\la_1=\la_2= 1$, and  $A \neq \Id$, while $f (x)=Ax+ae_1$, where $e_1$ is the eigenvector for  $A$. In this case the matrix logarithm $B$ is similar to degenerate Jordan cell.  The lines $f^t(x), t\in\rr$ are the integral curves for the autonomous system $\dot x=B x+b e_1$. Since $Bx $ is a  real multiple of $ e_1 $, the right side of the equation \eqref{dyn2} is the multiple of  $e_1$, so these curves are straight lines parallel to $e_1$. The line 
$ L=\{ -b e_2 +t e_1,t\in\rr\}$ is the set of fixed points for $f$.\\

{\bf Type 5.} This is the case when  $\la_1=\la_2= 1$,    $A \neq \Id$, and $f(x)=Ax+ue_1+ve_2$, where $e_2$ is the root vector for $A$ and $v\neq 0$. In this case $f$ has no fixed points. One can see that the integral curves corresponding to $f$ are parabolas obtained from each other by  parallel translations: 

Notice that matrix logarithm of $A$ is equal to $B=A-E$ and  $B^2=0$. 

Therefore the system $\dot x= B x+\be$ with initial value $x(0)=x_0$, has the solution \beq x(t)=x_0 +(B x_0 +\be)t + B  \dfrac{t^2}{2} \eeq

Denoting  $ ue_1+ve_2=b $, we get $\be = (I -\dfrac{B}{2} ) \cdot b$ and
$x(t)=B b \dfrac{t^2}{2} + (b - \dfrac{1}{2} B b +B x_0 ) t +x_0 $, while the vector field for $f$ is
\beq
\label{bfor5}
\dot x = B x +b - \frac{1}{2} B b \mbox{\qquad   or   \qquad}\dot x = (A-I) x + \left(\frac{3}{2}E -\frac{1}{2} A\right)b.
\eeq

Taking into account that for $x= \xi e_1 +\eta e_2$, $B x= \eta e_1$,  we see, that the right side in\eqref{bfor5}
$\eta e_1 +(u - v/2)e_1 + ve_2 $
does not depend on $\xi$ and vanishes $0$ if $v=0$ and  $\eta=-u$,  which corresponds to  Type 4. 

Therefore if $f$ belongs to the Type 5 the vector field has no stationary points and is preserved by  translations by  $te_1$,  so the minimal value of $\|\dot{x}\|$ is $|v| \cdot \|e_2\| \cdot |\sin \al_{12}|$, where $\al_{12}$ is the angle between
$e_1$ and $e_2$.

\begin{lem}\label{nofp}
Suppose that under the conditions of Theorem \ref{gen}, all the maps $f_n$ belong to the Type 1. Then there is such sequence of non-degenerate affine maps $h_n$ satisfying the conditions of Theorem \ref{gen} that their fixed points $y_n=\fix (h_n) \notin \overline{D} $.

\end{lem}

If $\ga$ is not a straight line segment, there is such a ball $B_1 \subset D$, that $\ga \cap \dot B_1 \neq 0$ and the set $\{ n: \fix f_n \IN \mC B_1\}$  is infinite.

By the condition (i) of the Theorem \ref{gen} there is such affine map  $g$, that $g(\ga) \subset \ga'$ and $g(B_1) \subset D$. Then, if $\fix f_n = x_n \in \mC B$, then $\fix(g^{-1} \cdot f \cdot g)= g^{-1} (x_n) \in \mC  g^{-1} (B_1) \subset \mC D$.

Thus  all the fixed points of the sequence of maps $f_n'=g^{-1} \cdot f_n \cdot g$ lie in the complement of $D$.

If $y=T  x + C$, then the fixed points  $y_n$ of the map $f_n'(x)$ are given by the equation $y_n = T^{-1} (x_n - C)$ and the map $f_n'$ is given by the equation $f_n'(x)= T^{-1} A_n T (x - y_n) + y_n$. 

At the same time the eigenvalues of the matrix   $A_n'$ are the same as the ones of $A_n$, and the sequence $f_n' \to \Id$.

Notice that for sufficiently large $n$ $f_n (g (a)) \IN g(\ga)$. Since $f_n$
 has no fixed points in $\ga$,  $f_n (g(\ga)) \cap g(\ga) = \ga(f_n(g(a_0), g(a_1)))$.
 
 Therefore $f_n'(\ga) \cap \ga = \ga((f_n'(a_0), a_1))$ and the sequence $f_n'$ satisfies the conditions of Theorem \ref{gen}.\vse \\

 {\bf Proof of Theorem \ref{gen}}
 
 Let $f_n$ be the sequence of maps  satisfying the conditions (i),(ii) of the  Theorem \ref{gen}. 
 
 Without loss of generality we may assume that all $f_k$ belong to  one and  the same of the Types 1-5.
 
 If all $f_k$ belong to the Type  2 or 4 then the set $P(x,a_0)$ lies on the segment $l_k=[a_0, f_k^{N_k} (a_0)]$, and the sequence  $l_k$ converges  to the segment $[a_0,a_1]$, therefore $\ga=[a_0,a_1]$.\\
 
 Thus  we need to prove the statement of the Theorem \ref{gen} for  the case when $f_n$ belong to Type 1,3 or 5. 
 
 If $f_n$ belong to  the Type 3 or 5, then the maps $f_n$ as well as their associated vector fields have no fixed points.

 If all  $f_n$ belong to the Type 1, Lemma \ref{nofp} allows us to assume that  fixed points of the maps $f_n$ lie outside of $D$.\\

 Let $L_k$ denote the set $\{ f_k^t (a_0), 0 \le t \le N_k\}$. Since $P(k,a_0)\IN L_k$  and  $\lim\limits_{k\to\8}P(k,a_0)=\ga$, we have $\ga \subset \overline{\lim \limits_{k \to \infty}} L_k $.

The sets $L_k$ are the subarcs of  integral curves of  linear dynamical systems  $ \dot x = B_k x + b_k$, and the endpoints of $L_k$ are $a_0$ and $f_k^{N_k} (a_0)$.

Let $m_k= \max \{ \| B_k x +b_k\|, x \in D\}$. If we replace the right sides  $B_k x + b_k$ of respective equations  \ref{bfor1},\ref{dyn2},\ref{bfor5}  by $B_k' x + b_k'$,  where $B_k' = B_k \setminus m_k$  and $b_k' = b_k \setminus m_k$, we obtain a sequence of linear dynamical systems in  $D$, which have no stationary points in  $D$, and whose integral curves are the same as the ones for the systems $\dot x = B_k x +b_k$. At the same time  
$\max \{ \| B_k' x +b_k'\|, x \in D\}$ is equal to 1 and  by convexity of the function $\| B_k' (x) + b_k'\|$, is assumed at some point  $x\in\dd D$. 

Denote $g_k(x) = B_k'+ b_k'$. The affine map $g_k$ sends  $D$ to some ellipse $g_k(D)\IN D$ which is tangent to $\dd D$  at some point and which does not contain  $0$. The sequence of maps $g_k$ satisfies the conditions of Arcela's theorem and one can find a subsequence $g_{n_k}$ which converges uniformly on $D$ to some affine function $g_0(x)$.

By continuous dependence of  solutions of differential equations on their right sides, the solutions of the differential equations $\dot x (t) = g_n (x), x(0)=a_0$ converge uniformly with all their derivatives to the solution of the equation $\dot x (t) = g_0 (x), x(0)=a_0$, and the integral curves $L_{n_k}$ converge to the curve $L_0$.  The curve $L_0$  belongs to the class  $C^2$ if $\|g_0(x)\|\neq 0$ so we need to control zero points of $g_0(x)$.

For that reason we consider the limit $g_0 (D)$ of the sequence of  ellipses $g_n(D)$.  

If $g_0 (D)$ is a non-degenerate ellipse, then since  $g_0 (D) = \lim g_{n_k} (D)$, and  $g_{n_k} (D) \notin 0 $,  $g_0 (D)$ can contain $0$ only on its boundary. Since $\ga \subset \dot D$, $g_0 (\ga) \notin 0$ in this case.

If $g_0 (D)$ -- is a line segment, for which $0$ is its inner point, then $g_0^{-1} (0)$  is a chord  $\La$ in the disc $D$. If  $\ga\IN\La$ then  $\ga$ is a line segment. Otherwise $\ga$ contains a subarc $\ga'$, which is disjoint from  $\La$. By the condition (i) we may assume that $\ga' = S(\ga)$ for some affine mapping $S$. The arc $\ga'$ is contained in the integral curve of the equation $ \dot x = g_0 (x)$, which starts at the point $S(a_0)$. Since $\|g_0(x)\|\neq 0$ on $\ga'$, it belongs to the class $C^2$. Therefore $\ga$ is twice differentiable.

By Theorem of C.Bandt and A.S.Kravchenko \cite[Theorem 3]{BK}, $\ga$ is a segment of a parabola or straight line.\vse

{\bf Proof of Theorem \ref{fin}.}

Let $f_n=S_{\bi_n}\yy S_{\bj_n}$ be the sequence  converging to $\Id$ for which $f_n(\ga)\cap\ga\neq\0$. Since $f_n$ is close to $\Id$, the maps $f_n$ 
and  $f_n\yy$ preserve the orientation on $\ga$. Notice that for self-affine arcs the condition (i) of Theorem \ref{gen} holds automatically. Therefore, following the argument of Lemma~\ref{nofp}, the sequence $f_n$ can be chosen in such a way that for any $n$, $\fix(f_n)\cap\ga=\0$. Then up to permutation of $\bi$ and $\bj$ we may suppose that for any $n$, $S_{\bi_n}(\ga)\cap S_{\bj_n}(\ga)=\ga(S_{\bj_n}(a_0), S_{\bi_n}(a_1))$. Therefore $f_n(\ga)\cap\ga=\ga(f_n(a_0),a_1)$ and we can apply Theorem~\ref{gen} to complete the proof.\vse

 \begin{dfn}
Let $\ga_1,\ga_2$ be Jordan arcs in Rd. We say that $\ga_1$ and $\ga_2$  have {\em proper intersection}
if the set  $\ga_1\cap\ga_2$ is a non-degenerate subarc in $\ga_1$ and $\ga_2$  and one of its endpoints is an endpoint of $\ga_1$ and
the other is an endpoint of $\ga_2$.
\end{dfn}

\begin{cor}\label{nine3}
Let $\eS$ be a system of non-degenerate contracting affine mappings with a Jordan attractor $\ga$.
 Let $A_\da(\ga)$ be the set of subarcs
$\al=h(\ga)\cap\ga$ such that  $|\al|\ge\da$, 
 $h$ is an affine map, and  the arcs $h(\ga)$ and $\ga$ have regular intersection.  If the set $A_\da(\ga)$
is infinite, then  $\ga$ is a segment of parabola.
\end{cor}

\section{ The partition to elementary subarcs.}

\begin{thm}\label{T4}
 Let $\eS=\{ S_1,...,S_m \}$ be a system of contractive affine maps in 
 $\rr^2$ with Jordan attractor  $\ga$.
 If $\ga$ is different  from a segment of a parabola or straight line, то there is
a multizipper $\eZ$ such that the arc $\ga$ is one  of the components of the attractor of $\eZ$.
\end{thm}

\dok  \  We suppose the system $\eS$ is irreducible.  Let us order the maps $S_1,...,S_m$ so that 
$\ga_i\cap\ga_j\neq\0$ if and only if  $|i-j|=1$, while $a_0\in\ga_1$ and $a_1\in\ga_m$. For two points $x,y\in\ga$ we write,
 that $x<y$, if $y\in \ga(x,a_1)$.  \\

First we construct such finite set $\eP\IN\ga$,
whose points $a_0=p_0<p_1<...<p_{N-1}<p_N=a_1$ define a partition of $\ga$ to subarcs $\da_i, i=1,...,N$, satisfying the conditions \\
1. For any $\da_i$   and any $k=1,...,m$ there is  $\da_j$ such that $S_k(\da_i)\IN \da_j$;\\
2. For any $k_1,k_2=1,...,m$ and for any $\da_{i_1}, \da_{i_2}$, $S_{k_1}(\dot\da_{i_1})$ and $S_{k_2}(\dot\da_{i_2})$ are either equal or disjoint. \\

Let $\eG$ be the set of all affine mappings $g$ such that the set $\ga\cap g(\ga)$ contains a connected component which is a subarc
$\ga_g\subset\ga$,
whose endpoints are the points $g(a_i)$ and $a_j$, $i,j\in\{0,1\}$. Let  $\eP$  be the set consisting of $a_0, a_1$ and of points $g(a_{i})$, where $g\in\eG$, $i=0,1$, and
$g(a_i)\in \ga_g\cap\dot\ga$. Let $\eP_i$ be the set of those $p\in \eP\cap\dot\ga$, 
which are the endpoints of subarcs $\ga_g$, that do not contain $a_{1-i}$. Thus, $\eP=\{a_0,a_1\}\cup\eP_0\cup\eP_1$.\\

Notice two properties of  $\eP$, which directly  follow from its definition:

{\bf b1}. Let $g$ be affine map of  $\rr^2$ for which
$g(\ga)\subset\ga$. Then $\eP\cap\dot g(\ga)\IN g(\eP)$.

{\bf b2}. Let $g_1,g_2$ be two affine maps such that 
$g_1(\ga),g_2(\ga)$ are subarcs of
  $\ga$, having proper intersection. Then the endpoint
 of the subarc $g_1(\ga)$, contained in  $g_2(\dot{\ga})$,
lies in  $g_2(\eP)$, and vice versa.\\

In the case when a Jordan arc $\ga$  is the attractor of a system
of contracting affine maps  $\eS$, the conditions  {\bf b1} and {\bf b2} imply the properties:\\

 {\bf c1}. For any $j\in I$,  $\eP\cap\dot\ga_j\IN S_j(\eP)$;

 {\bf c2}. For any $1\le j\le m-1$, $S_j(\{ a_0,a_1 \}\cap\dot{\ga}_{j+1}\IN  g_{j+1}(
\eP)$  and $S_{j+1}(\{ a_0,a_1 \}\cap\dot{\ga}_{j}\IN  g_{j}(
\eP)$\\

\begin{lem}\label{l1}
Let a Jordan arc $\ga\IN\rr^2$ with endpoints  $a_0,a_1$
be the attractor of irreducible system  $ \eS =\{S_1, ..., S_m\}$ of contracting
affine maps, and $\ga$ is not a segment of a parabola or a straight line. Then:\\
 {\bf d1}. The set of limit points of $\eP$ is contained in  $\{a_0,a_1  \}$.\\
  {\bf d2}. There are such neighbourhoods  $U_i$ of the points  $a_{i}$, where $i=0,1$, that   $P_{1-i}\cap U_i=\0$, and\\
   {\bf d3}. If for some $k\in\{1,m\}$ and some $i,j\in\{0,1\}$, $S_k(a_i)=a_j$, then  $S_k$ is a bijection of $U_i\cap
 \eP_i$ to $S_k(U_i)\cap \eP_j$.
\end{lem}

\dok   First we show that the set $\eP$ has no limit points in $\dot\ga$.  Suppose there is a $c\in\dot\ga\cap\bar\eP$.
Then for one of the endpoints of $\ga$, say, for $a_{0}$,
there is a sequence $g_n\in{\eG}$,
such that $g_n(a_{0})\to c$.
 It follows from Corollary \ref{nine3}, that $\ga$ is a segment of a parabola, which contradicts the assumptions of the Lemma, so
 {\bf d1} is true. The same argument shows that  $a_1$ cannot be  a limit point of $\eP_0$ and  $a_0$  cannot be a limit point of $\eP_1$. Therefore there are such neighbourhood $U_i$ of the points $a_i$, that $\eP_{1-i}\cap U_i=\0$. Moreover, we choose $U_0,U_1$ in  such a way that  $\ga\cap U_0\IN \ga_1$ and  $\ga\cap U_1\IN \ga_m$.\\
 
 To check {\bf d3}, consider first the case when $S_1(a_1)=a_0$. If $p\in \eP_0\cap U_0$ and $p=g(a_i)$, then
 $S_1\yy \circ g \in\eG$ and $S\yy_1(p)\in\eP_1\cap S\yy_1(U_0)$. Conversely, if $p\in \eP_1\cap U_1$, and $p=g(a_i)$,  then $S_1\circ g\in\eG$ and $S_1(p)\in\eP_0\cap S_1(U_1)$. Therefore  $S_1$ defines a bijection $\eP\cap U_0\cap S_1(U_1)$ to $\eP\cap U_1\cap S\yy_1(U_0)$. Enumerating all possibilities:\\ 1.$S_1(a_0)=a_0,\ S_m(a_1)=a_1$;\\ 2.$S_1(a_0)=a_0,\ S_m(a_1)=a_0$;\\ 3.$S_1(a_0)=a_1,\ S_m(a_1)=a_1$;\\ 4.$S_1(a_0)=a_1,\ S_m(a_1)=a_0$,\\ we find the desired pairs of  neighborhoods for  each of the cases. $\blacksquare$

\begin{lem}\label{l2}
The set  $\eP$  contains a finite subset $\eP'$, which also satisfies {\bf c1} and {\bf c2}.

\end{lem}

 \dok
For each of the points $S_k(a_i)\in\dot\ga$, where $k\in I$ and $i=0,1$ we denote by $w(k,i)$ the connected component of the set $\ga_k\mmm\eP$, which has $S_k(a_i)$ as its endpoint, whereas for $S_k(a_i)=a_j$ we put $w(k,i)=U_j$.  Let $W_i=\bigcap\limits_{k\in I}S_k\yy(w(k,i))\cap U_i$.\\
Let $\eP'=\{a_0,a_1\}\cup \eP\mmm(W_0\cup W_1)$.\\

The set  $\eP'$ is finite, so we denote its elements by $a_0=p_0<p_1<...<P_M=a_1$, and the subarcs $\ga(p_{k-1},p_k)$--- by $\da_k$.

For any $j\in I$,  $S_j(\eP)\IN S_j(W_0\cup W_1)\cup S_j(\eP')$. At the same time the definition of $\eP'$
implies that $S_j(W_0\cup W_1)\cup S_j(\eP')=S_j(\{a_0,a_1\})$.
Therefore  $\eP'\cap\ga_j\IN S_j(\eP')$. Thus the set $\eP'$
satisfies the condition {\bf c1}. The condition {\bf c2} directly follows from the definition of $\eP'$.
 $\blacksquare$

\begin{lem}\label{l3}
Each of the subarcs $\da_i,i=1,...,M$ and $\ga_i, i\in I$ is an  union of   subarcs $S_j(\da_k)$ for some $j\in I$ and some $k\in \{1,...,M\}$ whose interiors are disjoint.
\end{lem}

\dok The system $\eS$ is irreducible, therefore each subarc 
$\ga_{j}, 1<j<m$ intersects two adjacent subarcs 
$\ga_{j-1},\ga_{j+1}$, so that
$\ga_{j}\setminus(\ga_{j-1}\cup\ga_{j+1})\neq\0$. For any subarc
$\bar{\ga}_{j}=\ga_{j}\setminus(\dot\ga_{j-1}\cup\dot\ga_{j+1})$
its enpoints  by {\bf c2}  are contained in  
$S_{j}(\eP')$; let them be the points
$S_j(p_{k_j}),S_j(p_{K_j})$. The arc
$\bar{\ga}_{j}$ has unique representation
$\bigcup\limits_{i=k_j}^{K_j-1} S_{j}(\da_{i})$. For each of the subarcs  $\ga_{j}\cap\ga_{j+1 }$  there are exactly two partitions: first, to the subarcs 
  $S_{j}(\da_{i})$ and second, to the subarcs
 $S_{j+1}(\da_{i})$;  choose one of them.
 Taking the  union over all subarcs and renumerating all the points, we get the desired partition for the whole $\ga$. By the property {\bf c1}, the partition we obtained is at the same time a partition for  each of the subarcs $\da_{k}$. $\blacksquare$

{\bf Proof of the Theorem \ref{T4}.} Now we can construct a Jordan multizipper, for which the components of the attractor will be the subarcs $\da_{j}$. Each of the subarcs $\da_{j},j=1,\ldots M$ is
a finite union of consequent subarcs  $S_{i}(\da_{k})$, which form a partition of $\da_j$. Therefore we can create a graph $\w G$ whose vertices are $u_j=\da_{j}$ and an edge  $e_{ij}$ is directed from  $u_i$ to $u_j$ if
there is such $S_k$, that $S_k(U_j)\IN \da_i$. \vse

\end{document}